\documentclass[11pt,a4paper]{article}
\usepackage{mathrsfs}
\usepackage{amsfonts,amssymb,amsmath,indentfirst,amsthm}
\usepackage{color}

\newcommand{\C}{\mathbb{C}}

\def\q \m#1#2{{\raise1pt\hbox{$#1$}\kern-1pt\big/
               \kern-1pt\raise-1pt\hbox{$#2$}}}

\def\Z{{\mathbb Z}}

\def\C{{\mathbb C}}

\def\1{{\bf 1}}

\def\<{\langle}
\def\>{\rangle}

\def\O{{\cal O}}

\def\be{\begin{equation}}
\def\ed{\end{equation}}
 \def\bes{\begin{equation*}}
\def\eds{\end{equation*}}
\def\1{\textbf 1}
 \makeatletter

\newcommand{\Rmnum}[1]{\expandafter\@slowromancap\romannumeral #1@}
\makeatother
\newfam\msbfam
\font\twelmsb=msbm10 at 12pt 
\font\sevenmsb=msbm10 at 7pt \font\fivemsb=msbm10 at 5pt
\textfont\msbfam=\twelmsb
\scriptfont\msbfam=\sevenmsb \scriptscriptfont\msbfam=\fivemsb
\newtheorem{thm}{Theorem}[section]

\newtheorem{lem}[thm]{Lemma}
\newtheorem{prop}[thm]{Proposition}

\theoremstyle{remark}

\theoremstyle{definition}
\newtheorem{defn}[thm]{Definition}

\newcommand{\g}{{\mathbf g}}

\newcommand{\comment}[1]{}
\numberwithin{equation}{section}
\begin{document}
\title{  A Vertex Algebra Commutant for the $\beta\gamma$-System and Howe pairs}
\author{Yan-Jun Chu, Fang Huang, Zhu-Jun Zheng \thanks{Supported in part by NSFC with grant Number
10471034 \& 10971071
 and Provincial Foundation of Innovative
Scholars of Henan.}} \maketitle
\begin{center}
\begin{minipage}{5in}
{\bf  Abstract}: Analogue to commutants in the  theory of
associative algebras, one can construct a new subalgebra of vertex
algebra known as a vertex algebra commutant. In this paper, for the
adjoint representation $V$ of Lie algebra $sl(2,\C)$, we describe a
commutant  of $\beta\gamma$- System $S(V)$ by giving its generators,
moreover, we get a new Howe pair of vertex algebras.
\end{minipage}
\end{center}
\section{Introduction}
 Let $W$ be a vertex algebra, and $U$ be its subalgebra,
 one
 can construct a new subalgebra, which is known as the commutant $Com(U,W)$
 of $U$ in $W$. In fact, this construction is a generalization of the coset
 construction in conformal field theory due to Kac-Peterson \cite{2} and
 Goddard-Ken-Olive \cite{3}, and was introduced by I. Frenkel and Zhu in \cite{4}
 in mathematics. It is an analogue to the ordinary commutant construction
 in associative algebra theory.

 In order to
 describe vertex algebra commutant  $Com(U,W)$ more clearly, we
 expect to find its generator set and the corresponding
 relations. But we don't know
 whether the  commutant $Com(U, W)$ is generated finitely, and how
 to find its generators. These
 are non-trivial problems.
 It's obviously that $U\subset Com(Com(U, W),W)$. If $U=Com(Com(U,
 W),W)$, the pair $U$ and $(Com(U, W))$ is called a Howe pair(\cite{14}).
 As in associate algebra, it should have applications in vertex algebra theory.
 Since up to now, we know little about these new objects, it should be
 useful  to construct new more examples in the vertex algebra
 category.

 According to the paper \cite{5}, the notion of a commutative circle
 algebra is abstractly equivalent to the notion of a vertex algebra,
 and there are the related correspondence between two theories. We
 shall refer to a commutative circle algebra simply as a vertex
 algebra throughout this paper.

 Let $\g$ be a finite dimensional complex semisimple Lie algebra,
 $V$ be a finite dimensional complex
 $\g-$module via the Lie algebra homomorphism $\rho:\g\longrightarrow
 EndV$.  The representation $\rho$ induces a vertex algebra
 homomorphism $\widehat{\rho}$ from $\O(\g, B)$ to $S(V)$, where $B$
 is the bilinear form $B(u,v)= -Tr(\rho(u)\rho(v))$ and $S(V)$
 is the $\beta\gamma-$system which was introduced in
 \cite{1}). If $V$ admits a
 symmetric invariant bilinear form $B'$, there is a vertex algebra
 homomorphism
 $\widehat{\psi}:\O(sl(2,\C),-\frac{dimV}{8}K)\longrightarrow
 S(V)^{\Theta_{+}}$, denote by
 $\mathcal{A}=\widehat{\psi}(\O(sl(2,\C),-\frac{n}{8}K))$(cf. \S 2 in
 \cite{7} in details). In \cite{13}, Zhu introduced a functor which
 assigns each vertex algebra $W$ to an associative algebra
 $A(W)$(known as the Zhu algebra of $W$). From \cite{4}, we know that
 the associative algebra $A(\O(\g,B))$ is isomorphic to the
 universal enveloping algebra $U(\g)$. Let $D(V)$ be the Weyl algebra
 of polynomial differential operators of $V$. It is well known that
 $A(S(V))$ is isomorphic to $D(V)$.
 Just like the classical commutant theory of associative algebras,
 how to describe the commutant $S(V)^{\Theta_+}$ of
 $\Theta=\widehat{\rho}(\O(\g, B))$ in $S(V)$
 is an important problem.

 As we know, there are some results on the descriptions of the
 commutant $S(V)^{\Theta_+}$. In \cite{7,10}, B. Lian and Linshaw
 studied the vertex algebra and invariant theory, and reduced the problem of
 describing $S(V)^{\Theta_+}$ to a problem in
 commutative algebra. They single out a certain category
 $\mathfrak{R}$ of vertex algebras equipped with a $\Z_{\geq
 0}-$ filtration  such that the associated graded objects are
 $\Z_{\geq 0}-$graded $\partial-$ rings. All vertex algebras of the form $S(V)$
 and $\O(\g,B)$ belong to the category $ \mathfrak{R}$, and so are their
 subalgebras $\widehat{\rho}(\O(\g,B))$ and $S(V)^{\Theta_+}$. In particular, the
 assignment $W \longmapsto gr(W)$ is a functor from
 $\mathfrak{R}$ to the category $\mathcal{R}$ of $\Z_{\geq 0}-$ graded
 $\partial-$rings, and
 the main object of study $S(V)^{\Theta_+}$ is sent to the
 $\partial-$ ring $gr(S(V)^{\Theta_+})$. It's lucky that the
 reconstruction property of the category $\mathcal{R}$ tells us that
 if one can find a generator set of $gr(S(V)^{\Theta_+})$, then he
 can construct a generator set of $S(V)^{\Theta_+}$. However, describing
 generators of $gr(S(V))^{\Theta_+}$ is much easier than that of $gr(S(V)^{\Theta_+})$
 in the invariant theory. Moreover, there is a canonical injection
 $\Gamma:gr(S(V)^{\Theta_+})\longrightarrow gr(S(V))^{\Theta_+}$, if
 $\Gamma$ is surjective, the generator set of $gr(S(V))^{\Theta_+}$ can be regarded as the generator
set of $gr(S(V)^{\Theta_+})$ , hence, one need to find the generator
set of $gr(S(V))^{\Theta_+}$. As an example in \cite{10}, taking
$\g=V=sl(2,\C)$, Linshaw  studied the subalgebras
$gr(S_{\beta}(V)^{\Theta_+})$ and
 $gr(S_{\gamma}(V)^{\Theta_+})$ of $gr(S(V)^{\Theta_+})$, and gave
 a complete description of vertex algebras
 $S_{\beta}(V)^{\Theta_+}$ and $S_{\gamma}(V)^{\Theta_+}$. Moreover,
 he showed that $\mathcal{A}$ is isomorphic to the current
 algebra $\O(sl(2,\C),-\frac{3}{8}K)$. In terms of $V=\g=sl(2,\C)$, in  \cite{7}, the authors
 used the properties of Gr$\ddot{o}$bner bases to
 prove that  $S(V)^{\mathcal{A}_+}=\Theta$ and obtained a Howe pair
$(\Theta, S(V)^{\Theta_+})$ in $S(V)$, based on above approach to
describe $S(V)^{\Theta_+}$. About the case that $\g$
 is abelian Lie algebra acting diagonally on a vector space $V$, Linshaw found a finite set
 of generators for $S(V)^{\Theta_+}
 $, and showed that  $S(V)^{\Theta_+}$ is a simple vertex algebra
 and a member of Howe pair (\cite{9}). More generally, if $\g=sl(n,\C), so(n,\C), sp(2n,\C)$ and $V$ are the copies of
 standard representations of $\g$, the authors used tools from commutative algebra and algebraic
 geometry, in particular, the theory of jet schemes, to
 describe
 $S(V)^{\Theta_+}$ and give some Howe pairs in vertex algebra context(\cite{11}).

 In this paper, based on the theory of vertex algebras and $\partial-$ rings in \cite{5, 7, 8, 10, 11},
 we also study the case of $V=\g=sl(2,\C)$.
 Under the related properties of Hilbert
 series, we find all finite generators of invariant ring
 $gr(S(V))^{\Theta_+}$, and then we give the finite
 generator set of $S(V)^{\Theta_+}$ explicitly.
 Moreover, we get a new Howe pair $(\mathcal{A}, S(V)^{\mathcal{A}_+})$ in $S(V)$.
 Here is detailed outline of the paper. Firstly, the action of
$\g\otimes\C[t]$ on $gr(S(V))$ induced by the adjoint representation
$V$, forms the invariant
 subalgebra $gr(S(V))^{\g\otimes\C[t]}$. Denote by $P=gr(S(V))$, by the theorem 5.9
 in \cite{11}, we get the $\partial-$ ring $P^{\g\otimes\C[t]}$ is generated
 by $P_0^{\g}$. In particular, the
 finite generator set of ring $P_0^{\g}$ is also the
 generator set of  $P^{\g\otimes\C[t]}$ as a $\partial-$
 ring. Secondly, we find finite generators and the
 corresponding relations of invariant ring  $P_{0}^{\g}$ by
 solving the Hilbert series of $V\oplus V^{*}$.
 Since the embeddings
$gr(\mathcal{A})\hookrightarrow
 gr(S(V)^{\Theta_+})\hookrightarrow gr(S(V))^{\g\otimes\C[t]}$ are
 both surjections, by the
 reconstruction property of $\partial-$ rings, we give the finite
 generator set of $S(V)^{\Theta_+}$. Moreover, we also get a new
 Howe pair $(\mathcal{A}, S(V)^{\mathcal{A}_+})$.

\section{Vertex Algebras and Some Examples}
In this section, we give a summary of vertex algebras for this paper
can be read easily. Please refree to papers \cite{5, 7, 10, 11} for
the details. We shall use such vertex algebra notions throughout
this paper.

 Let $V$
 be a vector space over $\C$, and
 $z, w$ be the formal variables. Denote the space of all linear maps
 $V\rightarrow V((z))$ by $EndV((z))$, where
 $$V((z)):=\left\{\sum\limits_{n\in \Z}v(n)z^{-n-1}\mid v(n)\in V,
 v(n)=0~\mbox{for}~n\gg 0\right\}.$$ Then each $v\in EndV((z))$ can be
 uniquely expressed as a formal series
 $v=v(z):=\sum\limits_{n\in \Z}v(n)z^{-n-1}$, where $v(n)\in EndV$.

 On the space $EndV((z))$, for $n\in \Z$, one can define n-th circle products
 as follows: For $u,v\in
 EndV((z))$, n-th circle products is defined by
 $$u(w)\circ_nv(w)=Res_{z=0}(u(z)v(w)\ell_{|z|>|w|}(z-w)^n-v(w)u(v)\ell_{|w|>|z|}(z-w)^n).$$
Here $\ell_{|z|>|w|}f(z,w)\in \C[[z,z^{-1},w,w^{-1}]]$ denotes the
power series expansion of a rational function $f(z,w)$ in the region
$|z|>|w|$.

The non-negative circle products are connected through the operator
product expansion (OPE) formula. For $u,v\in EndV((z))$, there are
$$ u(z)v(w)=\sum\limits_{n\geq
0}u(w)\circ_nv(w)(z-w)^{-n-1}+:u(z)v(w):,$$ or it is written as
$ u(z)v(w)\sim \sum\limits_{n\geq 0}u(w)\circ_nv(w)(z-w)^{-n-1},$
where $\sim$ means equal modulo the term $:u(z)v(w):$.
Here, $:u(w)v(w):$ is a well defined element of $EndV((z))$, called
the Wick product of $u$ and $v$, and there is
$:u(w)v(w):=u\circ_{-1}v.$ The other circle products are related to
this by $n!u(z)_{-n-1}v(z)=:\partial^{n}u(z)v(z):$ for non-negative
integers n, where $\partial$ denotes the formal differentiation
operator $\frac{d}{dz}$.
For $u\in EndV((z))$, there is
$1\circ_n u=\delta_{n,-1}u$ for all $n\in \Z$; $u\circ_n
1=\delta_{n,-1}u$ for $n\geq-1$.

A linear subspace $U\subset EndV((z))$ containing $1$ which is
closed under the circle products will be called a circle algebra. In
particular, $U$ is closed under $\partial$ since $\partial
u=u\circ_{-2}1$.
Let $U$ be a circle algebra, a subset $S=\{u_i|i\in I\}$ of $U$ is
called to generate $U$ if any element of $U$ can be written as a
linear combination of non-associative words in the letters $u_i,
\circ_n$ for $i\in I$ and $n\in \Z$. It is said that $S$ strongly
generates $U$ if any element of $U$ can be written as a linear
combination of words in the letters $u_i, \circ_n$ for $n<0$.
Equivalently, $U$ is spanned by the collection
$\{:\partial^{k_1}u_{i_1}(z)\partial^{k_2}u_{i_2}(z)\cdots\partial^{k_m}u_{i_m}(z):|
k_1,k_2,\cdots,k_m\geq 0\}$.

\begin{defn}
We say that $u,v\in EndV((z))$ circle commute if
$(z-w)^{N}[u(z),v(w)]=0$ for some $N\geq 0$. Here  $[\ ,\ ]$ denotes
the bracket. If $N$ can be choose to be zero. we say that $u,v$
commute. A circle algebra is said to be commutative if its elements
pairwise circle commute.
\end{defn}

The notion of a commutative circle algebra is abstractly equivalent
to the notion of a vertex algebra. For any commutative circle
algebra $U$, define
$$
\begin{array}{ll}
\pi:U\longrightarrow (End U)((t))\\
u\longmapsto \pi(u): \pi(u)v=\sum\limits_{n\in \Z}(u\circ_n
v)t^{-n-1},~\mbox{for}~ \forall v\in U.
\end{array}
$$
Then $\pi$ is an injective circle algebra homomorphism, and the
quadruple of structure $(U,\pi,1,\partial)$ is a vertex algebra.
Conversely, if $(V,Y,1,D)$ is a vertex algebra, the collection
$\{Y(v,z)|v\in V\}\subset EndV((z))$ is a commutative circle
algebra.

Next, to study the theory of commutants of vertex algebras, we
introduce to two examples of vertex algebras.

 \textbf{Example
1(current algebra)}. Let $\mathbf g$ be a Lie algebra equipped with
a symmetric invariant bilinear form $B$, and $\C[t,t^{-1}]$ be the
Laurent polynomial algebra over $\C$ with one indeterminate $t$,
affine Lie algebra
$\widehat{\mathbf g}=\mathbf g \otimes \C[t,t^{-1}]\oplus \C c$ with
bracket
$$
[u\otimes t^n,v\otimes t^m]=[u,v]\otimes
t^{m+n}+nB(u,v)\delta_{m+n,0}c,$$ where $c$ is the center of
$\widehat{\mathbf g}$.

Set $deg(u\otimes t^n)=n, deg(K)=0$, then $\widehat{\mathbf g}$ is
equipped with the $\Z-$ grading structure. Let $\widehat{\mathbf
g}_{\geq 0}\subset \widehat{\mathbf g}$ be the subalgebra of
elements of non-negative degree, and let $N(\mathbf
g,B)=\mathfrak{U}(\widehat{\mathbf g})\otimes_{\widehat{\mathbf
g}_{\geq 0}}\C$ be the induced $\widehat{\mathbf g}-$module, where
$\C$ is the 1-dimensional $\widehat {\mathbf g}_{\geq 0}-$module on
which $\mathbf g \otimes\C[t]$ acts by zero and $c$ acts by 1.
 For each $u\in \mathbf g$, let $u(n)$
be the linear operator on $N(\mathbf g,B)$ representing $u\otimes
t^n$, and put $u(z)=\sum\limits_{n\in \Z}u(n)z^{-n-1}\in
End(N(\mathbf g,B))((z))$. The collection $\{u(z)|u\in \mathbf g\}$
generates a vertex algebra in $End(N(\mathbf g,B))((z))$, which we
denote by $ \O(\mathbf g, B)$(\cite {4,5,7,10}). For any $u,v\in
\g$, the vertex operators $u(z),v(z)\in End(N(\mathbf g,B))((z))$
satisfy OPE \be u(z)v(w)\sim B(u,v)(z-w)^{-2}+[u,v](w)(z-w)^{-1}.
\ed Denote $\1$ by the vacuum vector $1\otimes 1\in N(\mathbf g,B)$,
then there is the conclusion
\begin{lem}
The creation map $\chi: \O(\g,B)\longrightarrow N(\g,B)$ sending $
u(z)\longmapsto u(-1)\1$ is an isomorphism of $\O(\g,B)-$modules.
\end{lem}

If $\g$ is simple, for $\lambda\in \C, \lambda\neq -\frac{1}{2}$,
$\O(\g,B)$ has a conformal element
$L_{\O}(z)=\frac{1}{2\lambda+1}\sum\limits_{i=1}^{dim\g}:u_{i}(z)u_{i}(z):$,
where $\{u_i|i=1,2,\cdots, dim\g\}$ is an orthonormal basis of $\g$
with respect to the killing form $K$.

\textbf{Example 2($\beta\gamma-$system)}. Let $V$ be a finite
dimensional vector space. Regarding $V\oplus V^*$as an abelian  Lie
algebra,
affine Lie algebra $\eta(V)=(V\oplus V^*)\otimes \C[t,t^{-1}]\oplus
\C \tau$ with bracket
\be [(x,x')\otimes t^n,(y,y')\otimes t^m]=(\langle
y',x\rangle-\langle x',y\rangle)\delta_{m+n,0}\tau,\ed for $x,y\in
V$, and $x',y'\in V^*$. Let $\sigma \subset \eta(V)$ be the
subalgebra generated by $\tau, (x,0)\otimes t^n$, and $(0,x')\otimes
t^m$ for $n\geq 0,m>0$. Let $\C$ be the 1-dimensional
$\sigma-$module on which $(x,0)\otimes t^n$ and $(0,x')\otimes t^m$
act trivially and the central element $\tau$ acts by the identity.
Denote the linear operators representing $(x,0)\otimes
t^n,(0,x')\otimes t^n$ on the induced module
$\mathfrak{U}(\eta(V))\otimes_{\mathfrak{U}(\sigma)}\C$ by
$\beta^{x}(n),\gamma^{x'}(n-1)$, respectively, for $n\in \Z$. The
power series \be \beta^{x}(z)=\sum\limits_{n\in
\Z}\beta^{x}(n)z^{-n-1},\gamma^{x'}(z)=\sum\limits_{n\in
\Z}\gamma^{x}(n)z^{-n-1}\ed generate a vertex algebra $S(V)$ in
$End(\mathfrak{U}(\eta(V))\otimes_{\mathfrak{U}(\sigma)}\C)((z))$,
called $\beta\gamma-$ system(\cite{1}), and the generators satisfy
OPE
\be
\beta^{x}(z)\gamma^{x'}(w)\sim <x',x>(z-w)^{-1},
\beta^{x}(z)\beta^{y}(w)\sim 0,\gamma^{x'}(z)\gamma^{y'}(w)\sim 0.
\ed

Suppose that $V$ is a n-dimensional $\g-$module via
$\rho:\g\longrightarrow EndV$, where $\g$ is a finite dimensional
Lie algebra. There is the following relation between these two
vertex algebras.
\begin{lem}
\label{2.3}(\cite{7}) The above map $\rho$ induces a vertex algebra
homomorphism $\widehat{\rho}:\O(\g,B)\longrightarrow S(V)$, where
$B$ is the symmetric invariant bilinear form
$B(u,v)=-Tr(\rho(u)\rho(v))$ for $u,v\in \g$.
\end{lem}
Here, let$\{x_1,x_2,\cdots,x_n\}$ be a basis of $V$ and
$\{x'_1,x'_2,\cdots,x'_n\}$ be dual basis of $V^*$, the vertex
algebra homomorphism $\widehat{\rho}$ send $u(z)$ to
\be\widehat{u}(z)=-\sum\limits_{i=1}^{n}:\beta^{\rho(u)(x_i)}(z)\gamma^{x'_i}(z):,
\forall u\in \g.\ed  And there is the following operator product
expansions \be \widehat{u}(z)\widehat{v}(w)\sim
B(u,v)(z-w)^{-2}+\widehat{[u,v]}(w)(z-w)^{-1},\ed for $u,v\in \g$.
Incidently,  $S(V)$ has a conformal element $L_{S}(z)=
\sum\limits_{i=1}^{dimV}:\beta^{x_i}(z)\partial\gamma^{x'_i}(z):$.

Analogous to the commutant construction in the theory of associative
algebra, one can has a way to construct vertex subalgebras of a
vertex algebra, known as commutant subalgebras.
\begin{defn}
Let $W$ be a vertex algebra and $U$ be any subset of $W$. The
commutant of $U$ in $W$, denote by $Com(U,W)$, is defined to be the
set of vertex operators $v(w)\in W$ commuting strictly with each
element of $U$, that is, $[u(z),v(w)]=0$ for all $u(z)\in U$.
\end{defn}
In terms of above circle product, $[u(z),v(w)]=0$ if and only if
$u(z)\circ_n v(w)=0$ for all $n\geq 0$, so there is also
\be\label{f2.6} Com(U, W)= \{v(z)\in W|u(z)\circ_n v(z)=0, \forall
u(z)\in U, n\geq 0\} \ed

Regard $W^{U_+}$ as the subalgebra of invariants in $W$ under the
action of $U$. If $\Theta$ is a vertex algebra homomorphic image of
a current algebra $\O(\g,B)$, $W^{\Theta_+}$ is just the invariant
space $W^{\g\otimes \C[t]}$.  According to the formula (\ref{f2.6}),
the action of $\Theta$ on $W$ is induced by the non-negative circle
products, hence we can write $Com(\Theta, W)$ as $W^{\Theta_+}$.

According to the vertex algebra homomorphism $\widehat{\rho}$ in
Lemma \ref{2.3}, there is the following definition
\begin{defn}
Let $\Theta$ be the subalgebra $\widehat{\rho}(\O(\g,B))\subset
S(V)$. The commutant algebra
$S(V)^{\Theta_+}=Com(\widehat{\rho}(\O(\g,B)),S(V))$ is called the
algebra of invariant chiral differential operators on $V$.
\end{defn}

For a simple Lie algebra $\g$ and $B(u,v)=\lambda K(u,v)$ for
$\lambda \in \C, \lambda\neq -\frac{1}{2}$, $S(V)^{\Theta_+}$ has
the conformal elements $L(z)=L_{S}(z)-\widehat{\rho}(L_{\O}(z))$.

\begin{lem}
\label{l2.6}(\cite{7}) The conformal weight-zero subspace
$S(V)^{\Theta_+}_{0}\subset S(V)^{\Theta_+}$ coincides with the
classical invariant ring $Sym(V^{*})^{\g}$.
\end{lem}

Let $V$ be a n-dimensional vector space, $D(V)$ be the Weyl algebra
of polynomial differential operators of $V$, then $D(V)$ has
generators $\beta^{x},\gamma^{x'}$, which are linear in $x\in V,
x'\in V^{*}$, and satisfies the commutation relations $[\beta^{x},
\gamma^{x'}]=<x', x>$.
If $V$ is a n-dimensional $\g-$module via $\rho:\g\longrightarrow
EndV$, there is an induced action $\rho^{*}$ of $\g$ on $D(V)$,
moreover, $\g$ acts on $D(V)$ by derivations of degree 0,
and we have $gr(D(V)^{\g})=gr(D(V))^{\g}=Sym(V\oplus V^{*})^{\g}$.
The action of $\g$ on $D(V)$ can be realized by inner derivations.
We have a Lie algebra homomorphism $\tau: \g\longrightarrow D(V)$
given in chosen basis by \be
\tau=-\sum\limits_{i=1}^{n}\beta^{\rho(u)(x_i)}\gamma^{x'_i}, \ed
which can be extended to a map $\mathfrak{U}(\g)\longrightarrow
D(V)$, and the action of $u\in \g$ on $ v\in D(V)$ is given by
$\rho^{*}(v)=[\tau(u),v]$. Thus
$D(V)^{\g}=Com(\tau(\mathfrak{U}(\g),D(V))$.

Let $V$ be a n-dimensional $\g$- module equipped with a symmetric
invariant bilinear form $B'$. If $\{x_1,x_2,\cdots,x_n\}$ is an
orthonormal basis of $V$ with respect to $B'$, and
$\{x'_1,x'_2,\cdots,x'_n\}$ is the corresponding dual basis of
$V^{*}$, there are the following two lemmmas in \cite{7,10}
\begin{lem}
There is a Lie algebra homomorphism $\psi:sl(2,\C)\longrightarrow
D(V)^{\g}$ given in an orthonormal basis with respect to $B'$ by the
formulas \be h\longmapsto
\sum\limits_{i=1}^{n}\beta^{x_i}\gamma^{x'_i}, e\longmapsto
\frac{1}{2}\sum\limits_{i=1}^{n}\gamma^{x'_i}\gamma^{x'_i},
f\longmapsto
-\frac{1}{2}\sum\limits_{i=1}^{n}\beta^{x_i}\beta^{x_i}, \ed where
 $\{e,f,h\}$ denote the standard generators of $sl(2,\C)$,
 satisfying
 $$
 [e,f]=h, [h,e]=2e, [h,f]=-2f.
 $$
\end{lem}
\begin{lem}
\label{l2.8} The homomorphism $\psi:sl(2,\C)\longrightarrow
D(V)^{\g}$ induces a vertex algebra homomorphism
$\widehat{\psi}:\O(sl(2,\C),-\frac{n}{8}K)\longrightarrow
S(V)^{\Theta_{+}}$ by \be h(z)\longmapsto
v^{h}(z)=\sum\limits_{i=1}^{n}:\beta^{x_i}(z)\gamma^{x'_i}(z): \ed
\be e(z)\longmapsto
v^{e}(e)=\frac{1}{2}\sum\limits_{i=1}^{n}:\gamma^{x'_i}(z)\gamma^{x'_i}(z):
\ed \be f(z)\longmapsto
v^{f}(z)=-\frac{1}{2}\sum\limits_{i=1}^{n}:\beta^{x_i}(z)\beta^{x_i}(z):
\ed where $K$ is the Killing form of $sl(2,\C)$.
\end{lem}

\section{Category $\mathfrak{R}$ and $\partial-$ Rings}
In this section, we introduce a certain category of vertex algebras
and a category of $\partial-$ rings. Here we refer to the related
theory of \cite{7,8,10,11,12}.

 Let $\mathfrak{R}$ be the category of pairs $(W, deg)$, where $W$ is
 a vertex algebra equipped with a $\Z_{\geq 0}-$filtration
 $$W_{0}\subset W_{1}\subset W_{2}\subset\cdots, W=\bigcup\limits_{k\geq
 0}W_{k},$$such that $W_{0}=\C$, and for $a\in W_{k}, b\in W_{l}$,
 there are
 \be \label{3.1}
 a\circ_n b\in W_{k+l},~\mbox{for}~ n<0,
 \ed
 \be \label{3.2}
 a\circ_n b\in W_{k+l-1},~\mbox{for}~ n\geq 0.
 \ed
Here $W_{k}=0$ for $k<0$. A non-zero element $a(z)\in W$ is said to
have degree $d$ if $d$ is the minimal integer for which $a(z)\in
W_{d}$. Morphisms in $\mathfrak{R}$ are morphisms of vertex algebras
which preserve the above filtration. Filtration on vertex algebras
satisfying (\ref{3.1}),(\ref{3.2}) is induced by Haisheng Li in
\cite{8} known as good increasing filtration. If $W$ possesses such
a filtration, it follows that the associated graded object
$gr(W)=\bigoplus\limits_{k>0}W_{k}/W_{k-1}$ is a $\Z_{\geq 0}-$
graded associative, commutative algebra with a unit $1$ under a
product induced by the Wick product on $W$. Moreover, $gr(W)$ has a
derivation $\partial$ of degree zero 0, and for each $a\in W_{d}$
and $n\geq 0$, operators $a\circ_n$ on $W$ induces a derivation of
degree $d-1$ on $gr(W)$. For each $d\geq 1$, we have the projection
$\phi_d:W_{d}\longrightarrow W_{d}/W_{d-1}\in gr(W)$.

If $u,v\in gr(W)$ are homogeneous of degree $r,s$, respectively, and
$u(z)\in W_r,v(z)\in W_s$ are vertex operators such that
$\phi_{r}(u(z))=u$ and $ \phi_s(v(z))=v$, it follows that
$\phi_{r+s}(:u(z)v(z):)=uv$.

Let $\mathcal{R}$ denote the category of  $\Z_{\geq 0}-$graded
commutative rings equipped with a derivation $\partial$ of degree
zero, which is called $\partial-$rings.

The prominent feature of $\mathcal{R}$ is that  vertex algebra $W\in
\mathcal{R}$ has the following reconstruction property. We can write
down a set of strong generators for $W$ as a vertex algebra just by
studying the ring structure of $gr(W)$. We say that the collection
$\{a_i|i\in I\}$ generates $gr(W)$ as a $\partial-$ring if the
collection $\{\partial^{k}a_i|i\in I, k\geq 0\}$ generates $gr(W)$
as a grading ring.

\begin{lem}(\cite{7})
\label{3.1}
 Let $W$ be a vertex algebra in $\mathcal{R}$. Suppose that $gr(W)$
 is generated as a $\partial-$ ring by a collection $\{a_i|i\in
 I\}$, where $a_i$ is homogeneous of degree $d_i$, choose vertex
 operators $a_i(z)\in W_{d_i}$ such that $\phi_{d_i}(a_i(z))=a_i$,
 then $W$ is strongly generated by the collection $\{a_i(z)|i\in
 I\}$.
\end{lem}

Define $k=k(W,deg)=Sup\{j\geq 1|W_{r}\circ_{n}W_{s}\subset
W_{r+s-j}, \forall r,s,n\geq 0\}$, it follows that $k$ is finite if
and only if $W$ is not abelian(cf. \cite{7}). For two vertex
algebras $\O(\g,B)$ and $S(V)$, there are $k(\O(\g,B),deg)=1,
k(S(V),deg)=2$.

\begin{lem}(\cite{7})
\label{3.2} Let $(W,deg)\in \mathfrak{R}$, and $k=k(W,deg)$ be as
above. For each $u(z)\in W$ of degree $d$ and $n\geq 0$, the
operator $u(z)\circ_{n}$ on $W$ induces a homogeneous derivation
$u(n)_{Der}$ on $gr(W)$ of degree $d-k$, defined on homogeneous
elements $v$ of degree $r$ by \be
u(n)_{Der}(v)=\phi_{r+d-k}(u(z)\circ_{n}v(z)). \ed Here $v(z)\in W$
is any vertex operator of degree $r$ such that $\phi_{r}(v(z))=v$.
\end{lem}


Let $(W,deg)\in \mathfrak{R}$, $\mathcal{C}$ be a subalgebra of  $W$
which is a homomorphic image of a current algebra $\O(\g,B)$.
Suppose that for each $u\in \g$, $u(z)\in \mathcal{C}$ has degree
$k$, so that the derivations $\{u(n)_{Der}|n\geq 0\}$ on $gr(W)$ are
homogeneous of degree 0. Then there is

\begin{lem}(\cite{7})
\label{l3.3} The derivations $\{u(n)_{Der}|n\geq 0\}$ form  a
representation of $\g\otimes \C[t]$ on $gr(W)$. Moreover, the
actions of $\g\otimes \C[t]$ on $W$ and $gr(W)$ are compatible in
the sense that for any $w(z)\in W$ of degree $r$, there are
$u(n)_{Der}\phi_{r}(w(z))=\phi_{r}\circ u(n)(w(z))$.
\end{lem}

Consider the invariant ring $gr(W)^{\mathcal{C}_{+}}$, since
$gr(W)^{\mathcal{C}_{+}}$ is closed under $\partial$,
 then $W^{\mathcal{C}_{+}}\hookrightarrow
W$ induces an injective ring homomorphism
$gr(W^{\mathcal{C}_{+}})\hookrightarrow gr(W)$ whose image clearly
lies in $gr(W)^{\mathcal{C}_{+}}$. So we have a canonical injective
homomorphism $\Gamma:gr(W^{\mathcal{C}_{+}})\hookrightarrow
gr(W)^{\mathcal{C}_{+}}$ as $\partial-$rings.


Let $W=S(V), \Theta=\widehat{\rho}(\O(\g, B))$, where $\g$ is
semisimple and $V$ is a finite dimensional $\g$-module. In the case,
$deg(\widehat{u}(z))=2=k$, so each $\widehat{u}(n)_{Der}$ is
homogeneous of degree 0 and $gr(S(V)$ is a $\g\otimes \C[t]-$module
by Lemma \ref{l3.3}. Denote by $P=gr(S(V))$, and we write the images
of $\partial^{k}\beta^{x}(z), \partial^{k}\gamma^{x'}(z)$ in $P$ as
$\beta_{k}^{x}$ and $\gamma^{x'}_{k}$, respectively. The action of
$\widehat{u}(n)_{Der}$ on the generators of $P$ is given by \be
\widehat{u}(n)_{Der}(\beta_{k}^{x})=C_{k}^{n}\beta_{k-n}^{\rho(u)(x)};
\widehat{u}(n)_{Der}(\gamma_{k}^{x'})=C_{k-n}^{n}\gamma^{\rho^{*}(u)(x')},
\ed where $C_{k}^{n}=k(k-1)\cdots (k-n+1)$ for $n, k\geq 0$,
$C_{k}^{0}=1, C_{k}^{n}=0$ for $n>k$.


Next, we shall state a conclusion in \cite{11}, which plays an
important role for this paper.

Let $V$ be a linear representation of $G$ (connected, reductive
linear algebraic group over $\C$ with Lie algebra $\g$). Choose a
basis $\{x_1,x_2,\cdots, x_n\}$ for $V^{*}$, so the regular function
ring $\O(V)\cong \C[x_1,x_2,\cdots, x_n]$,
and there is $P=gr(S(V))=\C[\beta_{k}^{x}, \gamma_{k}^{x'}|x\in V,
x'\in V^{*}, k\geq 0]$.

\begin{lem}(\cite{11})
\label{l3.4}
 Let $G$ be a connected, reductive algebraic group, and
let $V$ be a linear $G-$ representation such that $\O(V\oplus
V^{*})$ contains no invariant lines. Then $P^{\g\otimes \C[t]}$ is
generated by
$$
P_{0}^{G}=\C[\beta_{0}^{x}, \gamma_{0}^{x'}|x\in V, x'\in
V^{*}]^{G}\cong \O(V\oplus V^{*})^{G}
$$
as a $\partial-$ ring. In particular, if $\{f_1,f_2,\cdots,f_n\}$
generate $P_{0}^{G}$ as a ring, then $\{f_1,f_2,\cdots,f_n\}$
generate $P^{\g\otimes \C[t]}$ as a $\partial-$ ring.
\end{lem}
At the same times, we have the well-known theorem(cf.\cite{12})
\begin{lem}
\label{3.6} (Hilbert Theorem) If $G$ is a connected, reductive
algebraic group, then the invariants ring of polynomials
$\C[x_1,x_2,\cdots,x_n]^{G}$ is finitely generated.
\end{lem}
For $\g=sl(2,\C)$, we take $V=sl(2,\C)$ the adjoint representation.
Since special linear group $SL(2)$ is a connected, linear reductive
algebraic group with the Lie algebra $sl(2,\C)$. By Lemma
\ref{l3.4}, $P^{sl(2,\C)\otimes \C[t]}$ is generated by
$P_{0}^{sl(2,\C)}=\C[\beta_{0}^{x}, \gamma_{0}^{x'}|x\in V, x'\in
V^{*}]^{sl(2,\C)}\cong \O(V\oplus V^{*})^{sl(2,\C)} $ as a
$\partial-$ring. Using Lemma \ref{3.6}, we know that
$P_{0}^{sl(2,\C)}$
is finitely generated. Then $P^{\Theta_+}=P^{sl(2,\C)\otimes \C[t]}$
is finitely generated as $\partial-$ring.
To describe the generators of $S(V)^{\Theta_+}$, we need to describe
the generators of $P^{\Theta_+}$ by the Hilbert series of
$sl(2,\C)\oplus sl(2,\C)^{*}$.
\section {The Hilbert series and The generators of
$P_{0}^{sl(2,\C)}$}

In this section, we shall calculate the Hilbert series of
$P_{0}^{sl(2,\C)}$, then give the generators of $P_{0}^{sl(2,\C)}$,
which are also the generators of $P^{\Theta_+}$ as $\partial-$ ring.
Here, we refer to the related definitions and results in \cite{12}.

\begin{defn}
Let $G$ be a subgroup of general linear group $GL(n)$,
$T=\C[x_1,x_2, \cdots, x_n]$ be the polynomial ring, $G$ has an
action on $T$, denoted by $T^G$ for the ring of invariants of $G$.
$T=\bigoplus\limits_{d\geq 0}T_d$, where $T_d\subset T$ is the
subspace of homogeneous polynomials of degree $d$, then
$T^G=\bigoplus\limits_{d\geq 0}T^G\cap T_d$. There is a formal power
series in an indeterminate $t$ \be P(t)=\sum\limits_{d\geq
0}dim(T^G\cap T_d)t^d\in \Z[[t]] \ed is called the Hilbert series of
the grading ring $T^G$.
\end{defn}
In terms of the Hilbert series, there are the two important
theorems.

\begin{lem}(\cite{12})
\label{4.2} If $T^G$ is generated by homogeneous polynomials
 $f_1,f_2,\cdots, f_r$ of degree $d_1,d_2,\cdots,d_r$, then the
 Hilbert series of $T^G$ is the power series expansion at $t=0$ of
 rational function
 \be
 P(t)=\frac{F(t)}{(1-t^{d_1})(1-t^{d_2})\cdots (1-t^{d_r})}
 \ed
 for some $F(t)\in \Z[t]$.
 \end{lem}

 Let $V$ be any $n-$ dimensional representation of $sl(2,\C)$,
 consider the induced action of $sl(2,\C)$ on the polynomial ring
 $T(V)=\C[x_1,x_2, \cdots, x_n]$ of functions on $V$. Let
 $a_1,a_2,\cdot,a_n\in \Z$ be the weights (not necessarily distinct)
 of $sl(2,\C)$ in the weight-space decomposition of $V$, the function
 \be
 P(q;t)=\frac{1}{(1-q^{a_1}t)(1-q^{a_2}t)\cdots
 (1-q^{a_n}t)}=det\left(I_{V}-t
 \left(\begin{array}{ll}
 q~~~ 0\\
 0~~ q^{-1}
 \end{array}\right)_{V}\right)^{-1}
 \ed
is called the $q-$Hilbert series of the representation $V$. There is
the relations between $q-$Hilbert series and Hilbert series as
follows

\begin{prop}(\cite{12})
\label{p1} The invariant ring $T(V)^{sl(2,\C)}$ has Hilbert series
\bes P(t)=Res_{q=0}(q-q^{-1})P(q;t). \eds Equivalently, if
$P(q;t)=\sum\limits_{m\in\Z}c_{(m)}(t)q^m$, then
$P(t)=c_{(0)}(t)-c_{(-2)}(t)$. \end{prop}

Next, we shall compute the Hilbert series of
$P_{0}^{sl(2,\C)}=\C[\beta_{0}^{e},\beta_{0}^{f},\beta_{0}^{h},
\gamma_{0}^{e'},\gamma_{0}^{f'},\gamma_{0}^{h'}]^{sl(2,\C)}$. Take
$V=sl(2,\C)$ is the adjoint representation of $sl(2,\C)$, then
$P_{0}^{sl(2,\C)}$ is the the polynomial ring of functions on
$sl(2,\C)\oplus sl(2,\C)^{*}$. Obviously, we know that
$\{2,2,0,0,-2,-2\}$ is the set of all weights in the weight-space
decomposition of the representation $sl(2,\C)\oplus sl(2,\C)^{*}$.
Associated to the definition of $q-$ Hilbert series,  we obtain the
$q-$ Hilbert series of the representation $sl(2,\C)\oplus
sl(2,\C)^{*}$ as follows \be P(q;t) =\frac{1}{(1-q^2t)^2 (1-t)^2
(1-q^{-2}t)^2} \ed

By Proposition \ref{p1} and some calculus, we can give the Hilbert
series of the representation $sl(2,\C)\oplus sl(2,\C)^{*}$.

\begin{prop}
\label{p2} For the representation $sl(2,\C)\oplus sl(2,\C)^{*}$, the
Hilbert series of the invariant ring
$P_{0}^{sl(2,\C)}=\C[\beta_{0}^{e},\beta_{0}^{f},\beta_{0}^{h},
\gamma_{0}^{e'},\gamma_{0}^{f'},\gamma_{0}^{h'}]^{sl(2,\C)}$ is \be
P(t)=\frac{1}{(1-t^2)^3}\ed
\end{prop}
Proof. According to the Proposition \ref{p1}, we need to write
$P(q;t)$ as a formal series $\sum\limits_{m\in\Z}c_{(m)}(t)q^m$.
Since the $q-$ Hilbert series
$$
\begin{array}{lll}
P(q;t) =\frac{1}{(1-q^2t)^2 (1-t)^2 (1-q^{-2}t)^2}
=\frac{1}{(1-t)^2} (\frac{q^4}{(1-q^2t)^2(q^2-t)^2})
\end{array}
$$
After calculus, $P(q;t)$ can be decomposed into the following form
$$
\begin{array}{ll}
P(q;t)=\frac{1}{(1-t)^2}\left(\frac{2t^2}{(1-t^2)^3}\frac{1}{1-q^2t}+\frac{1}{(1-t^2)^2}\frac{1}{(1-q^2t)^2}+
\frac{2t}{(1-t^2)^3}\frac{1}{q^2-t}+\frac{t^2}{(1-t^2)^2}\frac{1}{(q^2-t)^2}\right)
\end{array}
$$
Using the following expansions of rational functions
$$\frac{1}{1-u}=\sum\limits_{n=0}^{\infty}u^n;\frac{1}{1+u}=\sum\limits_{n=0}^{\infty}(-1)^n
u^n;\frac{1}{(1-u)^2}=\sum\limits_{n=0}^{\infty}(n+1)u^n,$$ we  get
the following expansion
$$
\begin{array}{ll}
P(q;t)=&\frac{1}{(1-t)^2}
(\frac{2t^2}{(1-t^2)^3}\sum\limits_{n=0}^{\infty}t^nq^{2n}+
\frac{1}{(1-t^2)^2}\sum\limits_{n=0}^{\infty}(n+1)t^nq^{2n}\\
&+\frac{2t}{(1-t^2)^3}\sum\limits_{n=0}^{\infty}t^nq^{-2n-2}+\frac{t^2}{(1-t^2)^2}\sum\limits_{n=0}^{\infty}(n+1)t^nq^{-2n-4})\\
&=\sum\limits_{m\in\Z}c_{(m)}(t)q^m
\end{array}
$$
Hence there are
$$\begin{array}{lll}
C_{(0)}(t)=\frac{1}{(1-t)^2}\left(\frac{2t^2}{(1-t^2)^3}+\frac{1}{(1-t^2)^2}\right)=\frac{1}{(1-t)^2}\frac{1+t^2}{(1-t^2)^3};
C_{(-2)}(t)=\frac{1}{(1-t)^2}\frac{2t}{(1-t^2)^3}.
\end{array}
$$
Finally, we obtain the Hilbert series
$$
\begin{array}{ll}
P(t)&=C_{(0)}(t)-C_{(-2)}(t)=\frac{1}{(1-t)^2}\left(\frac{1+t^2}{(1-t^2)^3}-\frac{2t}{(1-t^2)^3}\right)\\
&=\frac{1}{(1-t^2)^3}.
\end{array}
$$

By Lemma \ref{4.2}, and Proposition \ref{p2},
$P_{0}^{sl(2,\C)}=\C[\beta_{0}^{e},\beta_{0}^{f},\beta_{0}^{h},
\gamma_{0}^{e'},\gamma_{0}^{f'},\gamma_{0}^{h'}]^{sl(2,\C)}$ has
three generators as a ring and these generators are all homogeneous
elements of degree 2. and using Lemma \ref{l3.4}, these generators
are also generators of $P^{\Theta_+}$ as $\partial-$ring.

\begin{prop}
\label{p4.5} There are the following three elements of degree 2 \be
v^{e}=4(\gamma_{0}^{h'}\gamma_{0}^{h'}+\gamma_{0}^{e'}\gamma_{0}^{f'}),
v^{f}=
-\frac{1}{16}(\beta_{0}^{h}\beta_{0}^{h}+4\beta_{0}^{e}\beta_{0}^{f}),
v^{h}=\beta_{0}^{h}\gamma_{0}^{h'}+\beta_{0}^{e}\gamma_{0}^{e'}+\beta_{0}^{f}\gamma_{0}^{f'},\ed
belong to $P_{0}^{sl(2,\C)}$.
\end{prop}

Proof. According to the Lemma \ref{l2.8}, we know that $v^{e}(z),
v^{f}(z), v^{h}(z)\in S(V)^{\Theta_+}$ with the same degree 2. If
$V=\g=sl(2,\C)$, then we can get an orthonormal basis
$\{\frac{h}{2\sqrt{2}}, \frac{e+f}{2\sqrt{2}},
\frac{e-f}{2\sqrt{-2}}\}$ with respect to the Killing form $K$ of
$sl(2,\C)$, the corresponding dual basis is $\{2\sqrt{2}h',
\sqrt{2}(e'+f'), \sqrt{-2}(e'-f')\}$. Hence there are
$$
\begin{array}{lll}
v^{e}(z)&=\frac{1}{2}(:\gamma^{2\sqrt{2}h'}(z)\gamma^{2\sqrt{2}h'}(z):
+ :\gamma^{\sqrt{2}(e'+f')}(z)\gamma^{\sqrt{2}(e'+f')}(z):\\
&\hspace{0.4cm}+:\gamma^{\sqrt{-2}(e'-f')} \gamma^{\sqrt{-2}(e'-f')})\\
&=4(:\gamma^{h'}(z)\gamma^{h'}(z):+:\gamma^{e'}(z)\gamma^{f'}(z):),
\end{array}
$$
$$
\begin{array}{lll}
v^{f}(z)&=-\frac{1}{2}(:\beta^{\frac{h}{2\sqrt{2}}}(z)\beta^{\frac{h}{2\sqrt{2}}}(z):+
:\beta^{\frac{e+f}{2\sqrt{2}}}(z)\beta^{\frac{e+f}{2\sqrt{2}}}(z):\\
&\hspace{0.4cm}+:\beta^{\frac{e-f}{2\sqrt{-2}}}(z)\beta^{\frac{e-f}{2\sqrt{-2}}}(z):)\\
&=-\frac{1}{16}(:\beta^{h}(z)\beta^{h}(z):+4:\beta^{e}(z)\beta^{f}(z):)
\end{array}
$$
$$
\begin{array}{lll}
v^{h}(z)&=:\beta^{\frac{h}{2\sqrt{2}}}(z)\gamma^{2\sqrt{2}h'}(z)+:\beta^{\frac{e+f}{2\sqrt{2}}}(z)\gamma^{\sqrt{2}(e'+f')}(z):\\
&\hspace{0.4cm}+:\beta^{\frac{e-f}{2\sqrt{-2}}}(z)\gamma^{\sqrt{-2}(e'-f')}(z):\\
&=:\beta^{h}(z)\gamma^{h'}(z):+:\beta^{e}(z)\gamma^{e'}(z):+:\beta^{f}(z)\gamma^{f'}(z)
\end{array}
$$
Denote by $\mathcal{A}=\widehat{\psi}(\O(sl(2,\C),-\frac{3}{8}K))$.
Since $gr(\mathcal{A})\hookrightarrow
gr(S(V)^{\Theta_+})\hookrightarrow
gr(S(V))^{\Theta_+}=P^{\Theta_+}$, so $\phi_2(v^{u}(z))\in
P_{0}^{sl(2,\C)}$, where $u=e,f,h$, hence we know $v^{e}, v^{f},
v^{h}$ belong to $P_{0}^{sl(2,\C)}$.

As we expect, there is the following conclusion
\begin{prop}
\label{p4.6} three elements $
v^{e}=4(\gamma_{0}^{h'}\gamma_{0}^{h'}+\gamma_{0}^{e'}\gamma_{0}^{f'}),
v^{f}=
-\frac{1}{16}(\beta_{0}^{h}\beta_{0}^{h}+4\beta_{0}^{e}\beta_{0}^{f}),
v^{h}=\beta_{0}^{h}\gamma_{0}^{h'}+\beta_{0}^{e}\gamma_{0}^{e'}+\beta_{0}^{f}\gamma_{0}^{f'},$
are algebraically independent in $P_{0}^{sl(2,\C)}$.
\end{prop}
Proof. Assume that there is a non-zero polynomial $H(y_1,y_2,y_3)\in
\C[y_1,y_2,y_3]$ such that $H(v^{e}, v^{f}, v^{h})=0$. We can
restrict to the subspace: $\gamma_{0}^{e'}=\beta_{0}^{e}=0;
\gamma_{0}^{f'}=\beta^{0}_{f}$, then there are
$v^{e}=4\gamma_{0}^{h'}\gamma_{0}^{h'};
v^{f}=-\frac{1}{16}\beta_{0}^{h}\beta_{0}^{h};
v^{h}=\beta_{0}^{h}\gamma_{0}^{h'}+\beta_{0}^{f}\gamma_{0}^{f'}$,
then define a map
$$
\begin{array}{ll}
\hspace{3cm}\mu: \C^3\longrightarrow \C^3\\
(\beta_{0}^{f}, \beta_{0}^{h}, \gamma_{0}^{h'})\longmapsto (v^{e},
v^{f},v^{h})=(4\gamma_{0}^{h'}\gamma_{0}^{h'},-\frac{1}{16}\beta_{0}^{h}\beta_{0}^{h},
\beta_{0}^{h}\gamma_{0}^{h'}+\beta_{0}^{f}\gamma_{0}^{f'}).
\end{array}
$$
Since one can separate variables for the map $\mu$, so it is a
surjective, hence $v^{e}, v^{f},v^{h}$ don't satisfy any identity
relation, of course, there isn't any polynomial relation of $v^{e},
v^{f},v^{h}$, then there is a contradiction. Therefore, it dosen't
exist non-zero polynomial $H(y_1,y_2,y_3)$ such that $H(v^{e},
v^{f}, v^{h})=0$. Three elements $
v^{e}=4(\gamma_{0}^{h'}\gamma_{0}^{h'}+\gamma_{0}^{e'}\gamma_{0}^{f'}),
v^{f}=-\frac{1}{16}(
\beta_{0}^{h}\beta_{0}^{h}+4\beta_{0}^{e}\beta_{0}^{f}),
v^{h}=\beta_{0}^{h}\gamma_{0}^{h'}+\beta_{0}^{e}\gamma_{0}^{e'}+\beta_{0}^{f}\gamma_{0}^{f'}$
are algebraically independent in $P_{0}^{sl(2,\C)}$.
\section{ The Main Results about the Commutant $S(V)^{\Theta_+}$}
According to given propositions in section 4,  if $\g=sl(2,\C)=V$,
we obtain some results as follows.

For the commutant $S(V)^{\Theta_+}$, we know that the weight-zero
subspace $S(V)^{\Theta_+}_0=Sym(V^{*})^{sl(2,\C)}$ by Lemma
\ref{l2.6}.

\begin{prop}
The conformal weight-zero subspace $S(V)^{\Theta_+}_0$ is generated
only by the element
$:\gamma^{e'}(z)\gamma^{f'}(z):+:\gamma^{h'}(z)\gamma^{h'}(z):$.
\end{prop}

Proof. At first, we know
$S(V)^{\Theta_+}_0=Sym(V^{*})^{sl(2,\C)}\cong \C[\gamma^{e'},
\gamma^{f'}, \gamma^{h'}]^{sl(2,\C)}$. Since the Hilbert series of
$V=sl(2,\C)$ is $P(t)=\frac{1}{1-t^2}$ in \cite{12}, so the
invariant ring $\C[\gamma^{e'}, \gamma^{f'},
\gamma^{h'}]^{sl(2,\C)}$ is only generated by an element of degree
$2$. It is easy to check
$\gamma^{e'}\gamma^{f'}+\gamma^{h'}\gamma^{h'}$ is an element with
degree $2$ in $\C[\gamma^{e'}, \gamma^{f'},
\gamma^{h'}]^{sl(2,\C)}$. Thus $S(V)^{\Theta_+}_0$ is generated by
the element
$:\gamma^{e'}(z)\gamma^{f'}(z):+:\gamma^{h'}(z)\gamma^{h'}(z):$ as
an algebra.
\begin{thm}
\label{t1} The invariant ring $P_{0}^{sl(2,\C)}$ is generated by
three elements
$v^{e}=4(\gamma_{0}^{h'}\gamma_{0}^{h'}+\gamma_{0}^{e'}\gamma_{0}^{f'}),
v^{f}=
-\frac{1}{16}(\beta_{0}^{h}\beta_{0}^{h}+4\beta_{0}^{e}\beta_{0}^{f}),
v^{h}=\beta_{0}^{h}\gamma_{0}^{h'}+\beta_{0}^{e}\gamma_{0}^{e'}+\beta_{0}^{f}\gamma_{0}^{f'}$,
therefore the invariant ring $P^{\Theta_+}$ is generated by three
elements $v^{e}, v^{f}, v^{h}$ as $\partial-$ring.
\end{thm}
Proof. Using Lemma \ref{l3.4} and  Proposition \ref{p2}, \ref{p4.5},
 \ref{p4.6}, we can get above result immediately.

Let $\mathcal{A}$ be the image of the homomorphism $\widehat{\psi}$
in Lemma \ref{2.3}. i.e.
$\mathcal{A}=\widehat{\psi}(\O(sl(2,\C),-\frac{3}{8}K))$, we have
\begin{thm}
\label{t2} the commutant $S(V)^{\Theta_+}$ is generated strongly by
three elements $$ v^{e}(z)
=4(:\gamma^{h'}(z)\gamma^{h'}(z):+:\gamma^{e'}(z)\gamma^{f'}(z):),$$
$$
v^{f}(z)=-\frac{1}{16}(:\beta^{h}(z)\beta^{h}(z):+4:\beta^{e}(z)\beta^{f}(z):),$$
$$
v^{h}(z)=:\beta^{h}(z)\gamma^{h'}(z):+:\beta^{e}(z)\gamma^{e'}(z):+:\beta^{f}(z)\gamma^{f'}(z),
$$
and there is  $S(V)^{\Theta_+}=\mathcal{A}$.
\end{thm}
Proof.  Since $\phi_2(v^{e}(z))=v^{e};\phi_2(v^{f}(z))=v^{f};
\phi_2(v^{h}(z))=v^{h}$ and we have known
$v^{e}(z),v^{f}(z),v^{h}(z)\in S(V)^{\Theta_+}$, so
$v^{e},v^{f},v^{h}\in gr(\mathcal{A})\hookrightarrow
gr(S(V)^{\Theta_+})$. By Theorem \ref{t1}, $P^{\Theta_+}$ is
generated by three elements $v^{e}, v^{f}, v^{h}$ as
$\partial-$ring, hence $gr(S(V)^{\Theta_+})=
gr(S(V))^{\Theta_+}=P^{\Theta_+}$.  According to Lemma \ref{3.1}, we
know $S(V)^{\Theta_+}$ is generated strongly by three elements $
v^{e}(z)
=4(:\gamma^{h'}(z)\gamma^{h'}(z):+:\gamma^{e'}(z)\gamma^{f'}(z):),$
$
v^{f}(z)=-\frac{1}{16}(:\beta^{h}(z)\beta^{h}(z):+4:\beta^{e}(z)\beta^{f}(z):),$
$
v^{h}(z)=:\beta^{h}(z)\gamma^{h'}(z):+:\beta^{e}(z)\gamma^{e'}(z):+:\beta^{f}(z)\gamma^{f'}(z)$.

However,  $A$ is also generated by these three elements $v^{e}(z),
v^{f}(z), v^{h}(z)$ as a subalgebra of $S(V)^{\Theta_+}$, so
$S(V)^{\Theta_+}=\mathcal{A}$.
\begin{thm}
\label{t3} The commutant
$S(V)^{\mathcal{A}_+}=Com(\mathcal{A},S(V))$ and $\mathcal{A}$ form
a Howe pair in vertex algebra  context.
\end{thm}

Proof. Since $S(V)^{\Theta_+}=\mathcal{A}$, then
$\Theta=\widehat{\rho}(\O(sl(2,\C), -K)$ is a subalgebra of
$S(V)^{\mathcal{A}_+}$, thus there is
$$\mathcal{A}\subset Com(S(V)^{\mathcal{A}_+}, S(V))\subset
Com(\Theta, S(V))=S(V)^{\Theta_+}=\mathcal{A},$$ so we have
$\mathcal{A}=Com(S(V)^{\mathcal{A}_+}, S(V))$, i.e, $(\mathcal{A},
S(V)^{\mathcal{A}_+})$ forms a Howe pair in vertex algebra context.

\noindent Yan-Jun Chu, Zhu-Jun Zheng\\
 {\small Department of Mathematics\\
 South China University of
 Technology\\
 Guangzhou 510641, P. R. China \\
 and\\
 Institute of Mathematics\\
 Henan University\\  Kaifeng 475001, P. R.
 China\\
 E-mail: zhengzj@scut.edu.cn}
\\
 Fang Huang\\
 {\small Department of Mathematics\\
 South China University of
 Technology\\
 Guangzhou 510641, P. R. China}

\end{document}